\documentclass[12pt,a4]{amsart}
\usepackage{amssymb}
\usepackage{amsmath}
\usepackage{amscd}
\usepackage{graphicx}
\usepackage{subfigure}

\newtheorem{proposition}{Proposition}[section]
\newtheorem{theorem}[proposition]{Theorem}
\newtheorem{lem}[proposition]{Lemma}
\newtheorem{cor}[proposition]{Corollary}

\theoremstyle{remark}
\theoremstyle{definition}

\newtheorem{rem}[proposition]{Remark}

\numberwithin{equation}{section}

\newcommand{\Z} {\mathbb{Z}}
\newcommand{\N} {\mathbb{N}}
\newcommand{\R} {\mathbb{R}}
\newcommand{\C} {\mathbb{C}}

\newcommand{\X}{\mathcal{X}}

\renewcommand{\L}{\mathcal{L}}

\begin{document}

\title
[Invariant measure and Tarski's plank problem]
{Invariant measure of rotational beta expansion and Tarski's plank problem}
\author{Shigeki Akiyama \and Jonathan Caalim}
\email{akiyama@math.tsukuba.ac.jp \hspace{10mm} nathan.caalim@gmail.com
}
\address{Institute of Mathematics \&
Center for Integrated Research in Fundamental Science and Engineering, University of Tsukuba, 1-1-1 Tennodai, \\
Tsukuba, Ibaraki, Japan (zip:350-8571)}
\address{
Institute of Mathematics, University of the Philippines Diliman, 1101
Quezon City, Philippines
}
\maketitle

\begin{abstract}
We study invariant measures of a piecewise expanding
map in $\R^m$ defined by an
expanding similitude modulo lattice.
Using the result of Bang \cite{Bang:51} on the plank problem of Tarski,
we show that when the similarity ratio is not less than $m+1$, it has an
absolutely continuous invariant measure equivalent to the $m$-dimensional
Lebesgue measure, under some mild assumption on the fundamental domain.
Applying the method to the case $m=2$, we obtain an alternative proof of
the result in \cite{Akiyama-Caalim:16} together with some improvement.
\end{abstract}

\begin{section}{Introduction}

Let $1<\beta\in \R$.
Fix an isometry $M$ from the orthogonal group $O(m)$ of dimension $m$.
Let $\L$ be a lattice of $\R^m$
generated by $\eta_i\ (i=1,2,\dots,m)$ and choose a standard fundamental domain
$$
\X=\left\{\left. \xi+ \sum_{i=1}^m x_i \eta_i\ \right|\ x_i\in [0,1)\right\}
$$
with a fixed translation $\xi\in \R^m$. Then
$$
\R^m = \bigcup_{d\in \L} (\X + d)
$$
is a disjoint partition of $\R^m$.
Define a map $T:\X \rightarrow \X$ by $T(z)= \beta M z - d$
where $d=d(z)$ is the unique element in $\L$
satisfying $\beta M(z)\in \X +d$.
We obtain an expansion
\begin{eqnarray*}
z &=& \frac {M^{-1}(d_1)}{\beta} + \frac{M^{-1}(T(z))}{\beta}\\
  &=& \frac{M^{-1}(d_1)}{\beta}+ \frac{M^{-2}(d_2)}{\beta^2} + \frac{M^{-2}(T^2(z))}{\beta^2} \\
  &=& \sum_{i=1}^{\infty} \frac {M^{-i}(d_i)}{\beta^i}
\end{eqnarray*}
with $d_i=d(T^{i-1}(z))$.
We call $T$ the {\it rotational beta transformation} in $\R^m$. 

The rotational beta transformations extend the notion of beta expansions \cite{Renyi:57,Parry:60,Ito-Takahashi:74} and negative beta expansions \cite{Ito-Sadahiro:09, Kalle:14} to higher dimensions. The classical beta expansions are themselves generalizations of the decimal expansion. There is an enormous literature dealing with problems relating non-integer representations of numbers in $\R$ and in higher dimensional spaces. 
Canonical number systems \cite{Gilbert:82a, Katai-Szabo:75, Scheicher-Thuswaldner:01, Akiyama-Pethoe:02} give representations of higher dimensional points using positive integer bases with little redundancy, that is, the set of points with more than one representation has 
Lebesgue measure zero. Disregarding this redundancy, 
it is also of interest to study representability of all complex numbers closed to the origin 
using a finite set of complex digits. 
This theme is repeatedly studied,
for instance, in \cite{Komornik-Loreti:07,Safer:99}.
The notion of rotational beta expansions, 
in the case $m=2$, can be viewed as representations of complex numbers using a complex number base with digits in $\Z^2$,
with little redundancy. 
Among such expansions, we found a lot of systems whose associated symbolic dynamics is sofic 
in \cite{Akiyama-Caalim:16}, which immediately 
give rises to construction of self-similar tilings using a finite number of polygons and their translations.
Note that the construction of self-similar tilings is an interesting difficult problem 
on its own (\cite{Kenyon:96, Vince:99, Lagarias-Wang:96b, Kenyon-Solomyak:10}), which gives us a motivation of the study of rotational beta expansion. 

In the study of the ergodic properties of these number theoretic expansions, the most prior subject is the existence and the uniqueness of the absolutely continuous invariant probability measure (henceforth, referred as ACIM) of the underlying dynamical system. 
Renyi \cite{Renyi:57} proved that the beta transformation $T(x)=\beta x \mod 1$ admits a unique ACIM,  and it is equivalent to the Lebesgue measure. 
Later, Parry \cite{Parry:60} gave the explicit form of this invariant density. 
His idea is to find an explicit fixed point of the associated Perron-Frobenius operator. 
In the same line \cite{Parry:79, Parry:64, Gora:07, Ito-Sadahiro:09} 
gave the explicit forms of the ACIM of similar piecewise expanding maps. 
By the simple feature of these maps that the 
forward orbits of the discontinuities fall into a single point, 
it follows that the ACIM is unique and therefore the system 
becomes ergodic (c.f. \cite{Li-Yorke:78}). 
For negative beta expansion, the ACIM is not always 
equivalent to the Lebesgue measure. 
Liao-Steiner \cite{Liao-Steiner:12}  studied the gaps between 
the supports (intervals) of the ACIM. 
Considering an appropriate skew product of beta expansions, we can 
study randomized system and 
their invariant measures (see \cite{Dajani-deVries:07,Kempton:14}).

Summing up, we know a lot on one dimensional beta expansions. 
Ergodic study becomes drastically involved if we move to higher dimensions, 
by the effect of higher dimensional discontinuities. 
By this difficulty, 
there are not so many results on the basic properties of higher dimensional 
`handy' expansion, 
although there are many number theoretical trials 
mentioned above. 
This gives another motivation to  study 
rotational beta expansions. 
Since $T$ is a piecewise expanding map, a general theory in 
\cite{Keller:79, Keller:85, Gora-Boyarsky:89, Saussol:00, Tsujii:00, Buzzi-Keller:01, 
Tsujii:01} tells us that
there exists at least one and at most a finite number of 
ergodic ACIM's, whose supports have
disjoint interiors (c.f. \cite{Saussol:00}). 
On the other hand, 
uniqueness of ACIM fails to hold in general rotational beta expansions. 
In \cite{Akiyama-Caalim:16} we studied the case $m=2$ and $M\in SO(2)$
and showed among others that $T$ has an ACIM equivalent to $2$-dimensional
Lebesgue measure when $\beta$ is larger than a certain bound depending only
on $\X$.

A {\it strip} is the closed region sandwiched by two parallel hyperplanes
of co-dimension one.
Its {\it width} is the distance between the two hyperplanes.
Define the {\it width} of a non empty set $A\subset \R^m$
by the minimum width of strips which contain $A$, and denote it by $w(A)$.
Let $r(P)$ be the {\it covering radius} of a point set $P\subseteq \R^m$,
i.e., $r(P)=\sup_{x\in \R^m} \inf_{y\in P} \| x-y\|$.
The values $w(A)$ or $r(P)$ can be $\infty$ if no such value exists in $\R$.
Define a property 
\begin{itemize}
\item[(S)]For any $z\in \X$ there exists $n\in \N$ that
$2r(T^{-n}(z)+\L)\le \beta w(\X)$.
\end{itemize}
This property is confirmed in many cases in a stronger form that there exists
$n\in \N$ that for any $z\in \X$ we have $2r(T^{-n}(z)+\L)\le\beta w(\X)$.
For instance,  if $2r(\L)\le \beta w(\X)$, then the last condition is valid with $n=0$,
because $r(z+\L)=r(\L)$.  In this paper we shall prove

\begin{theorem}\label{m+1}
Assume that the property (S) holds. 
If $\beta\ge m+1$, then $(\X,T)$ has a unique ACIM. This ACIM is equivalent to the $m$-dimensional Lebesgue measure.
\end{theorem}

The proof of Theorem \ref{m+1} depends on a beautiful
result of Bang \cite{Bang:51},
which gave an affirmative answer to the conjecture of Tarski
on the covering of a convex body by strips.

\begin{rem}
\label{Width}
The property (S) guarantees that for any $z\in \X$,
$
\X \setminus T^{-n-1}(z)
$
can not contain an open ball of radius $w(\X)/2$.
In fact, from
$\beta M(T^{-n-1}(z)) = \beta M(\X) \cap V$ for 
$V:=\{ y+d\ | y\in T^{-n}(z), d\in \L\}$, 
if there is a ball of radius $\beta \omega(\X)/2$
in $\beta M(\X) \setminus V$ then its center is at least $\beta w(\X)/2$ away
from any point of $V$, giving a contradiction.
The property (S) is easy to show in many cases.
However (S) is not always true, since
it is easy to give a counter example when $T$ is not surjective for $m>1$.
\end{rem}

\begin{rem}
\label{Best}
For $m=1$ we have $r(\L)=w(\X)/2$ and the assumptions are reduced to $2\le \beta$.
On the other hand \cite{Akiyama-Scheicher:04} showed that the symmetric beta expansion
$$
U:x\mapsto \beta x -\left\lfloor \beta x +\frac 12 \right\rfloor
$$
in $[-1/2,1/2)$
does not have an ACIM equivalent to $1$-dimensional Lebesgue measure if $\beta<2$.
This shows that the bound $m+1$ in Theorem \ref{m+1}
is best possible for $m=1$.
\end{rem}

To prove Theorem \ref{m+1}, we introduce an idea to look at the {\it hole}
in $\X$ in dealing with ACIM, whose definition is found in \S \ref{Two}.
By the same idea, one can improve the lower bound when
an additional condition is satisfied.

\begin{theorem}\label{2}
Assume that the property (S) holds and there is an $\eta\in \L$ such that 
$$\beta M \X \subset \bigcup_{j\in \Z}(\X+j\eta).$$
If $\beta>2$, then $(\X,T)$ has a unique ACIM. 
This ACIM is equivalent to the $m$-dimensional Lebesgue measure.
\end{theorem}

As a corollary to this result, we discuss a family of 
rotational beta transformations defined in the complex plane,
which was introduced by Klaus Scheicher and Paul Surer \cite{Scheicher-Surer}.
Let $\beta>1$. Let $\zeta\in \C \setminus \R$ where $|\zeta|=\beta$. 
Consider $\X=\{x+(-\overline{\zeta})y|x,y\in [0,1)\}$. 
Let $T_{\zeta}$ be the rotational beta transformation on $\X$ given by
$z\mapsto \ \zeta z-d$ where $d=d(z)\in \L=\Z+(-\overline{\zeta})\Z$. 
For $\theta \in (0,\pi)$, define 
$$C(\theta):=\sqrt{\frac{1+\sqrt{1+4\sin^4(\theta)}}{2\sin^2 (\theta)}}.$$

\begin{cor} If $\beta>\max\{2,C(\theta)\}$, then $(\X, T_{\zeta})$ has a unique ACIM. 
This ACIM is equivalent to the 2-dimensional Lebesgue measure.
\end{cor}

One can check that, in particular, if 
$\theta \in [\sin^{-1}(2/\sqrt{15}),\pi-\sin^{-1}(2/\sqrt{15})]$
then the corollary holds for all $\beta>2$.

Specializing our technique to $m=2$, we obtain an alternative proof of Theorem 1 in
\cite{Akiyama-Caalim:16}.
Let $\eta_1,\eta_2, \xi\in \R^2$
and
$$
\X=\{ \xi+x_1 \eta_1+ x_2 \eta_2 \ |\ x_i, x_2 \in [0,1) \}$$ be a
fundamental domain of the lattice $\L=\Z \eta_1 + \Z \eta_2$ in $\R^2$.
We fix an isometry $M\in O(2)$ to define the map $T$.
Denote by $\theta(\X)\in (0,\pi)$ the angle between $\eta_1$ and $\eta_2$.
Define for $\theta=\theta(\X) \in (0,\pi/2]$,
$$
B_1=B_1(\theta):=
\begin{cases}
2 & \mbox{ if } \frac{1}{2} < \tan (\theta/2)\\
1 + \frac 2{1 + \sin(\theta/2)} & \mbox { if } \sin(\theta)<\sqrt{5}-2\\
\frac 32 + \frac {\cot^2(\theta/2)}{16}  +
\tan^2(\theta/2) & \mbox{ otherwise}
\end{cases}
$$
and
$$
B_2=B_2(\theta):=\begin{cases}
 1 + \frac {1}{\sin(\theta(\X)) \cos(\theta/2)} &
 \mbox{ if }  \frac{\pi}3<\theta\\
1 + \frac 2{1 + \sin(\theta/2)} & \mbox{ otherwise.}
\end{cases}
$$
For $\theta(\X) \in [\pi/2,\pi)$, define $B_i(\theta):=B_i(\pi-\theta)$, ($i=1,2$).

\begin{theorem}\label{B}
Assume that $m=2$ and (S) is valid.
If $\beta >B_1$ then $(\X,T)$ has a unique ACIM $\mu$.
Moreover, if $\beta>B_2$
then $\mu$ is equivalent to the 2-dimensional Lebesgue measure.
\end{theorem}

One can confirm the inequality $B_1\le B_2<3$ in Figure \ref{bound}.
\bigskip

\begin{figure}[htp]
	\centering
		\includegraphics[scale=1]{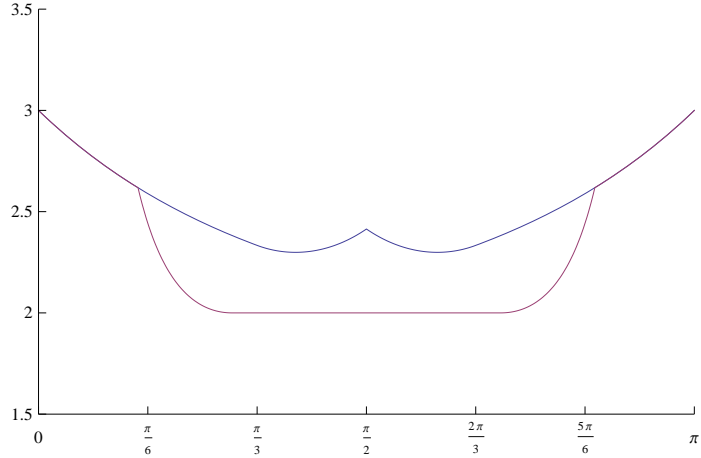}
		\caption{Comparison of $B_1$ and $B_2$\label{bound}}
\end{figure}

In fact, Theorem \ref{B} is an improvement of Theorem 1 in \cite{Akiyama-Caalim:16}.
First, it works not only for rotations but also for reflections. Second
the value of $B_1$ and $B_2$ are improved when $\min\{ \theta(\X), \pi-\theta(\X) \}$
is small, and became uniformly bounded by $3$, the bound given by
Theorem \ref{m+1}.
It is interesting that
we retrieve the bound in \cite{Akiyama-Caalim:16} when $\theta(\X)$ is close to $\pi/2$, though the methods are pretty different.

\begin{rem}
\label{Best2}
To illustrate Theorem \ref{B}, we discuss a simple family of rotational beta expansions.
Let $\X= [-1/2,1/2)\times [-1/2,1/2)$. 
Define $T$ to be the rotational beta expansion on $\X$ where 
the isometry $M$ is the identity. The map $T$ is a Cartesian product of
two identical maps $U:[-1/2, 1/2)\rightarrow [-1/2, 1/2)$ in Remark \ref{Best}.
If $\beta<2$, then the support of the ACIM of $U$ is contained in 
$Z=[-1/2,\beta/2-1]\cup [1-\beta/2,1/2]$. 
Clearly the support of the ACIM of $T$ is contained in $Z\times Z$ and the ACIM is
not  equivalent to the $2$-dimensional Lebesgue measure.
Let $\beta\le \sqrt{2}$ and 
consider the restrictions of $T$ in $Y_1$ and $Y_2$, where 
$Y_1$ is composed of eight squares of size $(\beta-1)/(\beta(\beta+1))$ and $(\beta-1)/2(\beta+1)$
and $Y_2$ is composed of eight $(\beta-1)/(\beta(\beta+1)) \times (\beta-1)/2(\beta+1)$
rectangles as in Figure \ref{TwoCompo}.
We easily confirm that in $Y_1$, a small square is mapped by $T$ into the diagonally opposite big square while a big square is mapped to the four corner squares. 
Hence, $T$ has an ACIM whose support is in $Y_1$. 
Likewise, the restriction of $T$ on $Y_2$ is well-defined and 
$T$ has at least two ergodic components. 
This is an interesting example that $1$-dimensional map
$U$ has a unique ACIM for every $\beta>1$ by the discussion of \cite{Li-Yorke:78}, 
but its Cartesian product $T$ has at least two
ACIM's when $\beta\le \sqrt{2}$. 
Consequently  the projections of different $2$-dimensional ACIM's of $T$ 
to the first coordinate are identical. 
This phenomenon is explained by 
the fact that $U$ is not totally ergodic;
by the ACIM of $U$, $U$ is ergodic but $U^2$ is not.
We can also construct isomorphic transformations where the 
fundamental domain is a rhombus of side length 1, $\theta(\X)\in (0,\pi)$ and 
whose diagonals meet at the origin. Meanwhile, we know that the unique ACIM preserved by $U$ 
is equivalent to the Lebesgue measure when $\beta \geq 2$ as mentioned in Remark \ref{Best}. 
Thus, when $\beta\ge 2$, $T$ admits a unique ACIM and this measure is equivalent to the 2-dimensional Lebesgue measure. Having these examples, we know  
$B_2\ge 2$ and $B_1>\sqrt{2}$ are necessary in Theorem \ref{B}. 
\begin{figure}[h]
\caption{$Y_1$ (Left) and $Y_2$ (Right)\label{TwoCompo}}
\centering
\includegraphics[width=0.7\textwidth]{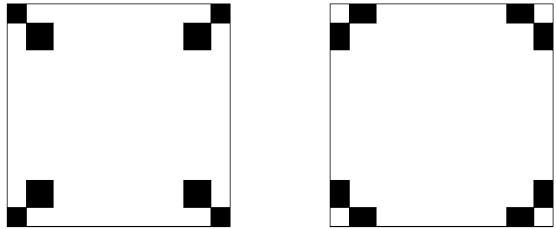}
\end{figure}
\end{rem}

\end{section}

\begin{section}[Two]{Proof of Theorem \ref{2}}

It is helpful for readers to start from the proof of Theorem \ref{2}.
In order to show that there exists an ACIM equivalent to Lebesgue measure, as
in \cite{Akiyama-Caalim:16}, it is enough to
prove that $\bigcup_{i=1}^{\infty} T^{-i}(z)$
is dense in $\X$ for any $z\in \X$.
Since $\X$ is bounded, the set of pairs $(x,r)\in \X
\times \R_{\ge 0}$ such that an open ball $B(x,r)$ is entirely contained in
$\X\setminus \bigcup_{i=1}^{n} T^{-i}(z)$
forms a compact set.  Hereafter such a ball is called a {\it hole of $n$-th level}.
Let $r_n$ be the maximum radius of the holes of $n$-the level. Fix $z\in \X$. Suppose that 
$\X\setminus T^{-n}(z)$ does not contain a ball of radius $w(\X)/2$ for some $n\in\N$.

Since
\begin{equation}
\label{Inv}
\beta M(T^{-n-1}(z)) = \beta M(\X) \cap \bigcup_{d\in \L} (T^{-n}(z) + d)
\end{equation}
and
$$
\beta M(\X) \subset \bigcup_{i\in \Z} (\X +i \eta),
$$
we consider the maximum open ball $B(x,r)$ contained in
$$
\bigcup_{i\in \Z} i\eta+ (\X \setminus T^{-n}(z)).
$$
By (\ref{Inv}), we have $r_{n+1}=r/\beta$.
Note that the largest hole $B(x,r)$ intersects
$\bigcup_{i\in \Z} \partial(\X)+i\eta$ at most once with
a hyperplane $(\X +i \eta) \cap (\X +(i+1)\eta)$.
In fact, if it intersects both
$(\X +i \eta) \cap (\X +(i+1)\eta)$ and
$(\X +(i+1) \eta) \cap (\X +(i+2)\eta)$,
then $B(x,r)$ contains a segment of length greater than $w(\X)$ and
the radius of the ball should be greater than $w(\X)/2$, 
which is not possible. We claim that
if $r>2 r_n$, then $B(x,r)\subset
\beta M(\X)$ must contain an element $y\in \bigcup_{i\in \Z}(T^{-n}(z)+i\eta)$.
If not, then since a hemisphere of $B(x,r)$
belongs to one side of the hyperplane,
there is a ball of radius not less than $r/2$ in $\X$ which does
not contain
an element of $T^{-n}(z)$ (see Figure \ref{Two}).
\begin{figure}[htp]
	\centering
		\includegraphics[scale=1]{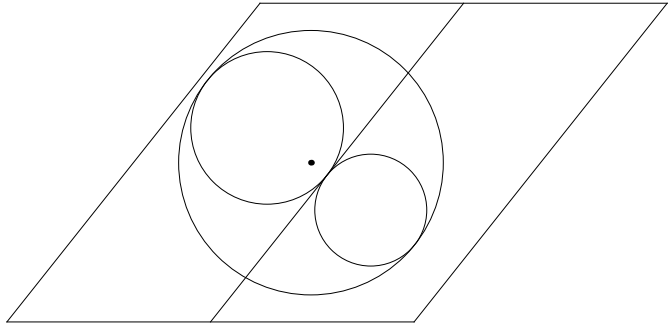}
\caption{The largest hole in $\beta M(\X)$ can not be doubled.\label{Two}}
\end{figure}
However $r/2>r_n$ gives a contradiction, showing the
claim.
Thus we have seen that $r\le 2 r_n$. From $r_{n+1}=r/\beta$, we get
$r_{n+1}\le 2r_n/\beta$. Iterating this procedure,
we see that $\lim_{n\rightarrow \infty}r_n=0$. This finishes the proof.
\end{section}

\begin{section}{Proof of Corollary}
Let $z\in \X$. Since $\beta^2>2$, 
the region $\zeta \X$ must contain at least two integer translates of $z$,
aligned with distance $1$. 
We find at least two points in $T^{-1}(z)$ aligned with distance $1/\beta$
in a line of slope being the argument of $-\overline{\zeta}$ of $\X$.
Let $L$ be the line segment 
connecting the midpoints of the longer sides of $\X$
and choose $x\in T^{-1}(z)$ closest to $L$.
Then, the segment of slope $-\overline{\zeta}$ 
starting from $x$ and ending at a point in $L$
is of length at most $1/(2\beta)$. 
Let us consider a circle of radius
$\beta w(\X)/2 =\beta \sin(\theta)/2$ inscribed in $\zeta \X$ as in Figure \ref{ZetaExp}. We consider the diameter parallel to the longer side of $\zeta \X$. Parallel to this diameter, we look at the chord of the circle extending to $x$ of length $2h$. Observing the 
right triangle of hypotenuse equal to the radius of the circle, a side of length $h$ 
and the other is of length at most $\sin(\theta)/(2\beta)$, we see
$$h \geq \sqrt{(\beta \sin(\theta)/2)^2-(\sin(\theta)/(2\beta))^2}.$$
Clearly, if $h>1/2$, then the chord is of length greater than 1 and it contains a translate $x+d$ with 
$d\in \Z$. Hence, an open ball of
radius $w(\X)/2$ in $\X$ must contain a point of ${T_{\zeta}}^{-2}(z)$.
Thus if $\beta>\max\{\sqrt{2},C(\theta)\}$ then $\X\setminus {T_{\zeta}}^{-2}(z)$ 
can not contain an open ball of radius $w(\X)/2$. 
Applying the proof of Theorem \ref{2} we get the desired result.

\begin{figure}[htp]
\includegraphics[scale=0.08]{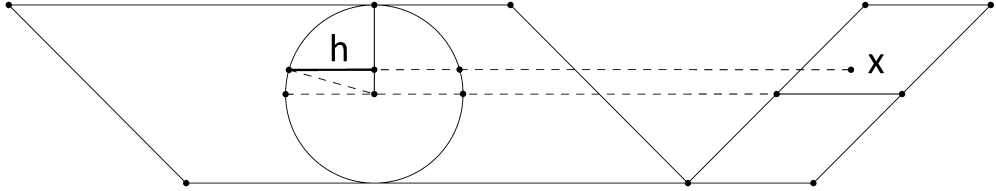}
\caption{$T_{\zeta}$ when $\zeta=\beta\exp(3\pi\sqrt{-1}/4)$\label{ZetaExp}}
\end{figure}

\end{section}

\begin{section}{Tarski's plank problem and dissection of a ball}

In this section we prepare a lemma to prove Theorem \ref{m+1}.
Motivated by the minimum dissection
of piecewise congruent polygons, A.~Tarski \cite{Tarski:32}
asked if a convex set $X$ is entirely covered by a set of strips (planks) of
width $w_i\ (i=1,\dots,k)$ whether $\sum_{i=1}^k w_i\ge w(X)$ must hold. This
problem is called Tarski's plank problem.
T.~Bang \cite{Bang:51} gave an affirmative answer to this question with an ingenious
proof.
Our lemma is closely related to this result.
Consider a unit ball $B$ in $\R^m$ and $k$ hyperplanes of co-dimension one which
cut the ball $B$ into convex cells. We are interested in the maximum of
the radius of open balls which can be inscribed in such cells.

\begin{lem}\label{Cut}
For arbitrary configurations of $k$ hyperplanes, there exists
an open ball of radius $1/(k+1)$ inscribed in a cell generated by the
$k$ hyperplanes and a $m$-dimensional unit ball.
The maximum radius of the balls in the cells
attains its minimum $1/(k+1)$ only when the
$k$ hyperplanes are parallel and totally
aligned by distance $2/(k+1)$.
\end{lem}

The first statement is found in Ball \cite{Ball:91} as well.

\begin{proof}
Consider a ball $W$ of radius $1-1/(k+1)=k/(k+1)$
and a unit ball $B$,
both centered at the origin.
For a hyperplane $H_i\ (i=1,2,\dots,k)$,
we associate a strip $P_i(\varepsilon)$ whose points are
of distance not greater than $1/(k+1)-\varepsilon$ from $H_i$
with a small
positive
constant $\varepsilon$. Since $$
k(2/(k+1)-2\varepsilon)<2k/(k+1)$$
by Bang's result, we have
$W \setminus \bigcup_{i=1}^k P_i(\varepsilon) \neq \emptyset$. Choosing
a center from this non-empty set, we find
an open ball of radius $1/(k+1)-\varepsilon$
entirely contained in a cell generated
by $B$ and $H_i\ (i=1,\dots,k)$.
Taking $\varepsilon$ smaller,
we obtain a decreasing sequence of non empty compact sets
$W \setminus \bigcup_{i=1}^k \mathrm{Int}(P_i(\varepsilon))$.
Thus the set $\bigcap_{\varepsilon>0}W \setminus \bigcup_{i=1}^k
\mathrm{Int}(P_i(\varepsilon))$ is non empty and
we obtain an open ball of radius $1/(k+1)$ entirely contained
 in a cell generated by $B$ and $H_i\ (i=1,\dots,k)$.

The remaining part of the proof depends more on the method by Bang.
A key of the proof of \cite{Bang:51} is that
$\bigcap (W \pm u_1 \pm u_2 \dots \pm u_k )\neq \emptyset$, where
$u_i$ is a vector  of length $1/(k+1)-\varepsilon$
perpendicular to the hyperplanes of the strip $P_i(\varepsilon)$.
The intersection is taken over all possible choices of signs.
This fact is valid even when $\varepsilon=0$, if not every $H_i$
are parallel.
In fact, without loss of generality, we may assume
that $H_i$ are parallel for $1\le i\le t<k$ and $H_j\ (j>t)$ are not parallel
to $H_1$. Then $\bigcap (W \pm u_1 \pm u_2 \dots \pm u_t )$ is simply
a lens shape $L:=(W - t u_1) \cap (W + t u_1)$ whose width
$2(k-t)/(k+1)$
is attained only
by the hyperplanes tangent to the two tops of the lens $L$ (see Figure \ref{Lens}).
\begin{figure}[htp]
	\centering
		\includegraphics[scale=0.7]{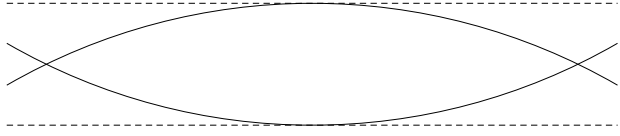}
\caption{Lens shape and its minimum strip\label{Lens}}
\end{figure}
From the proof of \cite{Bang:51}, we see that
$\bigcap (L \pm u_{t+1} \dots \pm u_k )\neq \emptyset$
contains a translation of $\kappa L$ with $\kappa>0$, because
and the minimal strip perpendicular to $u_j$ which contain $L$ has
width greater than $w(L)$ for $j>t$.
Following the proof of Bang,
we know that $P_i(0)\ (i=1,\dots, k)$ do not entirely cover $W$. Therefore
if not all $H_i$'s are parallel, we can find a ball of radius strictly larger than
$1/(k+1)$ entirely contained in a cell.
\end{proof}

\end{section}

\begin{section}{Proof of Theorem \ref{m+1}}
We proceed in the same manner as Theorem \ref{2}. We will prove
that $\bigcup_{i=1}^{\infty} T^{-n}(z)$ is dense in $\X$ for any $z\in \X$
by showing that the maximum radius $r_n$ of the holes of $n$-the level
converges to $0$.
From (\ref{Inv}),
we wish to get the maximum of the radius $r$ of the ball $B(x,r)
\subset \beta M(\X)$
which does not contain any element of
$\bigcup_{d\in \L} (T^{-n}(z) + d)$.
By Remark \ref{Width},
such a ball can not intersect two parallel hyperplanes which contain
$m-1$ dimensional faces of $\partial(\X)+d$ for $d\in \L$. Therefore
$B(x,r)$ is cut by at most $m$ hyperplanes $H_i\ (i=1,\dots,m)$, each of
which contain a face of some translation
$\partial(\X)+d$ with $d\in \L$.
Lemma \ref{Cut} shows there is
an inscribed ball of radius $r/(m+1)$ in some cell.
Thus we see within
$B(x,r)$, there is a ball of radius $r/(m+1)$ inscribed in $\X+d$ with $d\in \L$
having no point of $T^{-n}(z)+d$, which implies $r/(m+1)\le r_n$.
As $r_{n+1}=r/\beta$, we have
$r_{n+1}\le (m+1)r_n/\beta$, which gives the conclusion for $\beta>m+1$.

Finally we discuss the case $\beta=m+1$.
If $m\ge 2$, any pair of the $m$ hyperplanes
are not parallel and Lemma \ref{Cut} shows for any ball $B(x,r) \subset \beta M(\X)$
there is the largest hole of radius strictly greater than $r/(m+1)$, which
shows that the assertion is still valid for $\beta=m+1$.
For $m=1$, the problem is reduced to the study of the map
$$
T(x)=2x+a-\lfloor 2x+a \rfloor \text{ or }
T(x)=-2x+a-\lfloor -2x+a \rfloor
$$
defined in $[0,1)$ where $1$-dimensional Lebesgue measure is clearly the ACIM.
\end{section}

\begin{section}{Proof of Theorem \ref{B}}

Throughout this section, we assume $\xi=0$ because
the proof is identical for all $\xi$.
As in \cite{Akiyama-Caalim:16},
$B_i\ (i=1,2)$ are symmetric with respect to
$\theta(\X) \leftrightarrow \pi-\theta(\X)$,
we prove the case $0< \theta(\X)\le \pi/2$.
Let us show the bound $B_2$.
Similar to the proofs of Theorem \ref{2} and \ref{m+1},
in view of (\ref{Inv}) we wish to find the infimum radius $r$
that for every ball $B(x,r)\subset \beta M(\X)$, there exists a $d\in \L$
such that $B(x,r)\cap (\X+d)$
contains a ball of radius $r_n$ under the knowledge of $\theta(\X)$.
Then $r_{n+1}$ can be chosen to be $r/\beta$.
An equivalent dual problem is
to find the supremum of $r$ that one can find a ball $B(x,r)\subset \beta M(\X)$
that for all $d\in \L$, $B(x,r)\cap (\X+d)$ does not contain a ball of radius $r_n$.
Clearly the supremum is attained by the configuration that for at least one
$d\in \L$
the boundary of $B(x,r)\cap (\X+d)$ is circumscribed about a ball of radius $r_n$.
If there is only one line passing through $B(x,r)$, then the supremum
would be $2r_n$ as Theorem \ref{2}, or in Figure \ref{Two}. Since
$B_2>2$ for any $\theta(\X)$,
we will see that this case does not give the supremum.
By Remark \ref{Width}, the ball
$B(x,r)$ intersects at most two non parallel hyperplanes (lines).

Let us consider the case that $B(x,r)$ intersects the two lines. For example,
in Figure \ref{Perturbation}, the ball $B(0,r)$ is touching two
maximal balls of radius $r_n$ in the cells.
However, this configuration does not give the supremum, because
shifting a little the ball $B(0,r)$ along the bisector of the acute angle,
we easily see that every cell no longer contains a ball of radius $r_n$.
\begin{figure}[htp]
	\centering
		\includegraphics[scale=0.7]{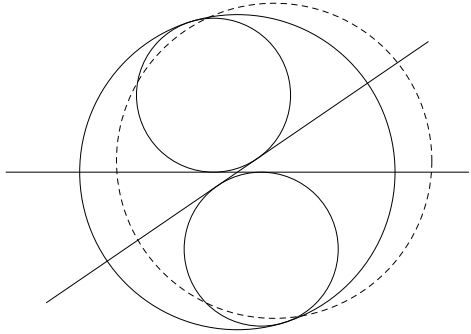}
\caption{Non supremum configuration\label{Perturbation}}
\end{figure}
Since a planar circle is determined by three intersections, by
similar small perturbation argument, in
the supremum configuration $\partial(B(x,r))$ must touch
$B(y_i,r_n)\ (i=1,2,3)$ where $B(y_i,r_n)$ are the balls inscribed in three distinct
$\X+d_i$ with $d_i\in \L$. Each $B(y_i,r_n)$ is a maximal ball within a
cell in $\partial(\X)+d_i$. There are exactly two possible supremum configurations
as depicted in Figure \ref{Small} and \ref{Large}.

\begin{figure}[ht]
\begin{center}
\subfigure[Small angle case\label{Small}]{
\includegraphics[scale=0.68]{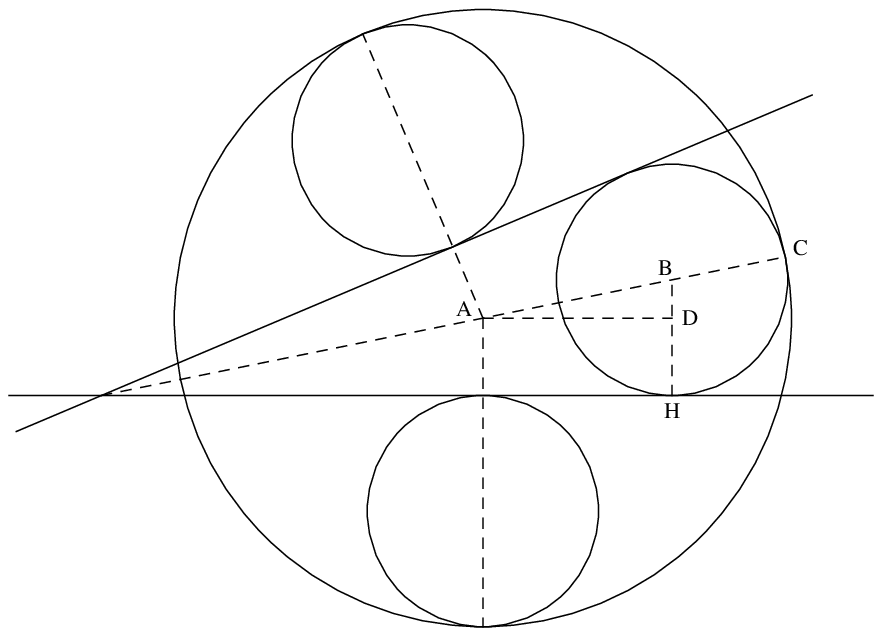}
}
\quad
\subfigure[Large angle case\label{Large}]{
\includegraphics[scale=0.68]{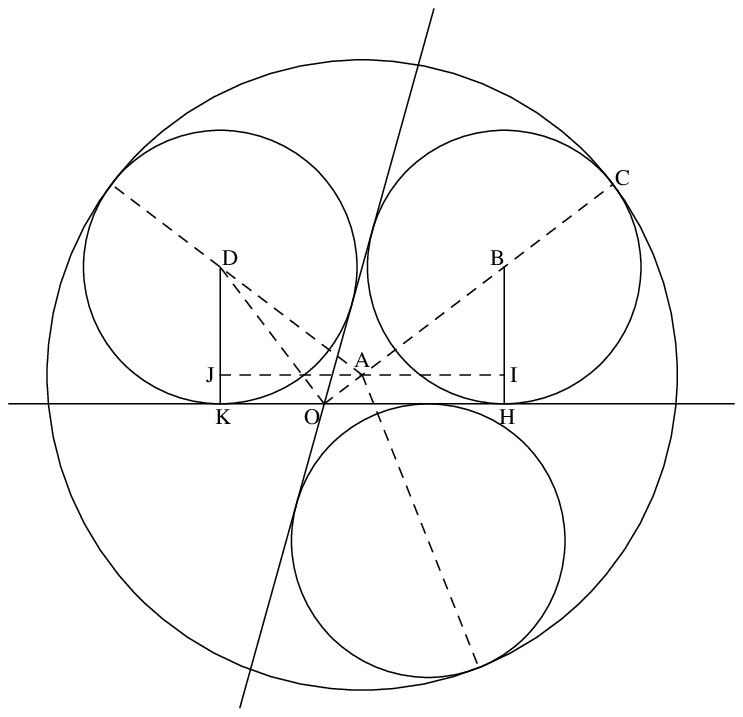}
}
\end{center}
\caption{Two possible supremum cases}
\end{figure}

Figure \ref{Small} is for $\theta(\X)<\pi/3$.
Points $A, B$ are the centers of the corresponding balls,
$BH$ is the perpendicular from $B$ to $x$-axis and
$AD$ is parallel to $x$-axis which intersects $BH$.
Since $AC=r$, $BC=BH=r_n$, putting $\ell=DH$, we have
$$
r=\ell+2r_n, \qquad \frac{r_n-\ell}{r_n+\ell}=\sin\left(\frac{\theta(\X)}2\right)
$$
and we obtain
$$
\frac {r}{r_n}= 1+\frac 2{1+\sin(\theta(\X)/2)}.
$$
For $\theta(\X)\ge \pi/3$ by the configuration of
Figure \ref{Small}, the above circle
intersects $x$-axis and the circle below intersects
the line of slope $\theta(\X)$.
Therefore we have to switch to Figure \ref{Large}.
The point $A$ is the center of the ball of radius $r$, whose radius is given by $AC=r$,
and the points $B,D$ are the ones of radius
$r_n$ and $O$ is the origin. The segments
$BH$, $DK$ are the perpendiculars from $B,D$
to $x$-axis. $JI$
is a segment passing through $A$ and parallel to $x$-axis
connecting the two perpendiculars. Then we have
\begin{eqnarray*}
KH &=& KO+OH = \frac {r_n}{\cot(\theta(\X)/2)} + \frac {r_n}{\tan(\theta(\X)/2)}\\
r&=&AB+r_n\\
AB \cos(\theta(\X)/2)&=&AI=AD \cos(\theta(\X)/2)=AJ= KH/2.
\end{eqnarray*}
Solving these equations we get
$$
\frac {r}{r_n}= 1 + \frac {\cot(\theta(\X)/2)+\tan(\theta(\X)/2)}{2 \cos(\theta(\X)/2)}
$$
which turned out to be equal to
$$
1+\frac{1}{\sin(\theta(\X))\cos (\theta(\X)/2)}
$$
appeared in Theorem 1.1 of \cite{Akiyama-Caalim:16}.
Next we prove the bound $B_1$. Our target is to show
the same statement in the proof of Theorem 1.1 in \cite{Akiyama-Caalim:16}:
\medskip

\begin{it}
If $\beta>B_1$, then for
any $\varepsilon>0$ there exists $z\in \X \setminus
\bigcup_{i=-\infty}^{\infty} T^i(\partial(\X))$ and
a positive integer $n$ that
$\bigcup_{i=0}^{n} T^{-i}(z)$ is an $\varepsilon$-covering of $\X$
\end{it}
\medskip

\noindent
by choosing a point $z$ very close to the origin.
In showing the bound $B_1$, it is practically the same as
assuming $z=0$ and the points of $\L$ can not be an inner point of $B(x,r)$.
Thus our geometry
problem has an additional restriction that the ball $B(x,r)$
can not contain the origin
as an inner point. Let us start with Figure \ref{Small}.
In this figure, the origin is outside of $B(x,r)$ and there is no change. However,
when we make $\theta(\X)$ gradually large, we find at some
angle the large circle passes through
the origin as in Figure \ref{Critical1}.
This happens when additionally
$$\frac {r_n}{\sin(\theta(\X)/2)}= r+\ell+r_n$$
holds, i.e., when $\sin(\theta(\X)/2)=\sqrt{5}-2$.
If
$\sin(\theta(\X)/2)>\sqrt{5}-2$, the supremum configuration is
switched to Figure \ref{TwoCircleMax}.


\begin{figure}[ht]
\begin{center}
\subfigure[First critical configuration\label{Critical1}]{
\includegraphics[scale=0.6]{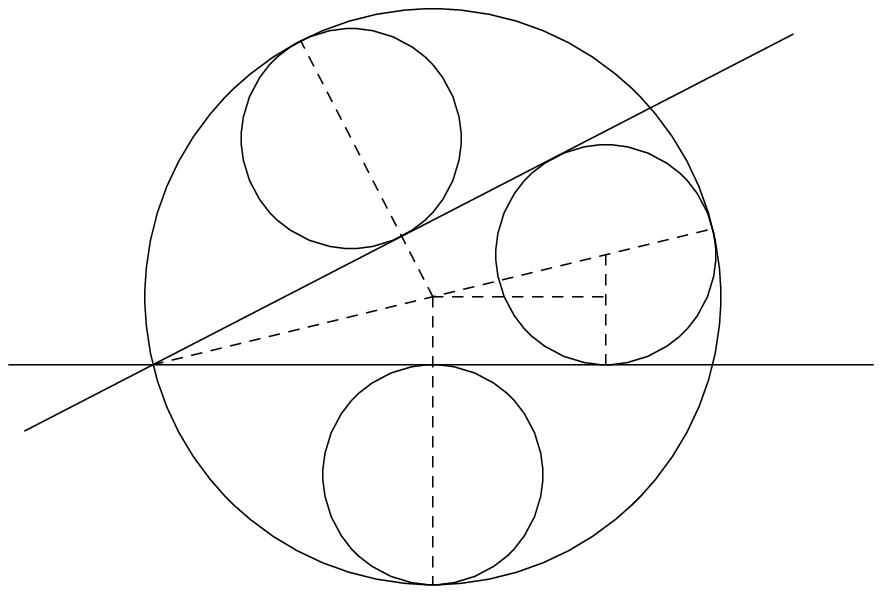}
}
\quad
\subfigure[Another configuration\label{TwoCircleMax}]{
\includegraphics[scale=0.6]{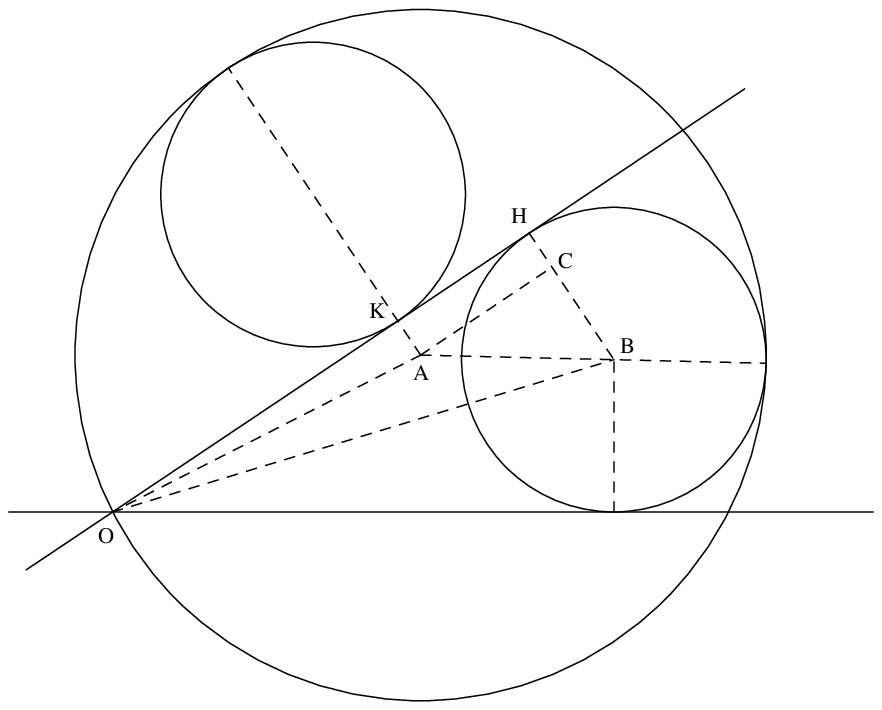}
}
\end{center}
\caption{Change of supremum configurations}
\end{figure}

This is the configuration that $\partial(B(x,r))$
passes through the origin and touches
$B(y_i,r_n)\ (i=1,2)$ where $B(y_i,r_n)$ are the balls inscribed in two distinct
$\X+d_i$ with $d_i\in \L$. Each $B(y_i,r_n)$ is a maximal ball in
a cell in $\partial(\X)+d_i$.
The point $A$ is the center of the ball of radius $r$ whose radius is given by $AO=r$,
$B$ is the center of the ball of radius $r_n$ and $O$ is the origin. The segments
$AK$ and $BH$
are perpendiculars from $A$ and $B$, respectively,
to the line of slope $\tan(\theta(\X))$.
The segment $AC$ is parallel to the line and intersects $BH$. Putting $AK=CH=\ell$, we have
\begin{eqnarray*}
r&=& \ell+2r_n\\
OK^2&=&r^2-\ell^2\\
AB&=&r_n+\ell\\
BC&=&r_n-\ell\\
AC&=&KH= \sqrt{4r_n \ell}\\
OH&=&OK+KH=\frac{r_n}{\tan(\theta(\X)/2)}.
\end{eqnarray*}
Summing up, we obtain
$$
\frac r{r_n}=
\frac 32 + \frac 1{16} \cot^2\left(\frac{\theta(\X)}2\right) +
\tan^2\left(\frac{\theta(\X)}2\right).
$$

This configuration works provided $OK\le OH$. The next critical configuration
is depicted in Figure \ref{Critical2}.

\begin{figure}[ht]
\begin{center}
\subfigure[Second critical configuration\label{Critical2}]{
\includegraphics[scale=0.6]{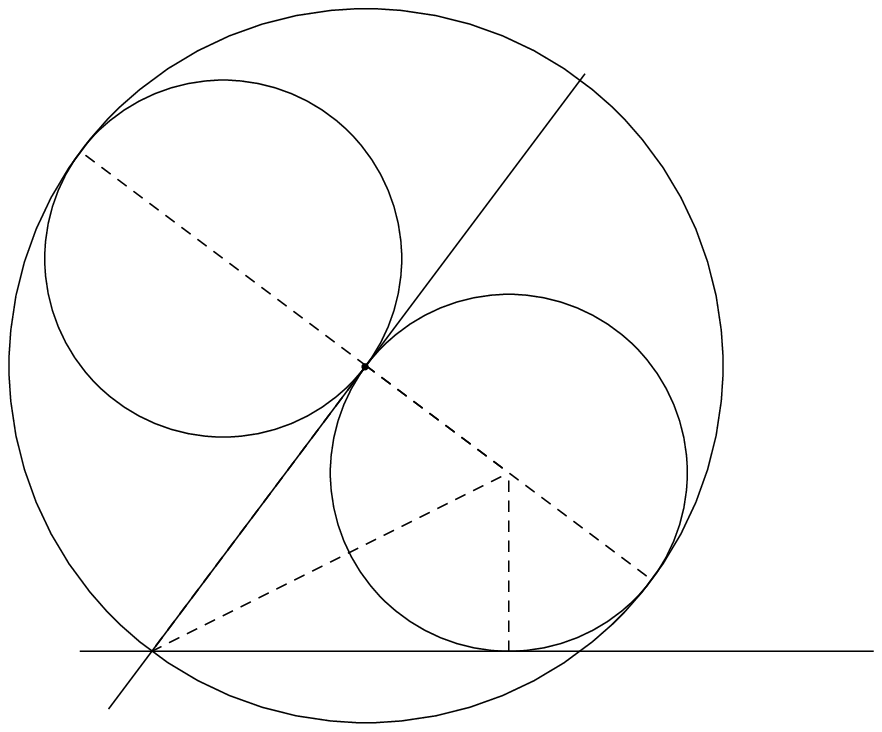}
}
\quad
\subfigure[Larger angle \label{Larger}]{
\includegraphics[scale=0.6]{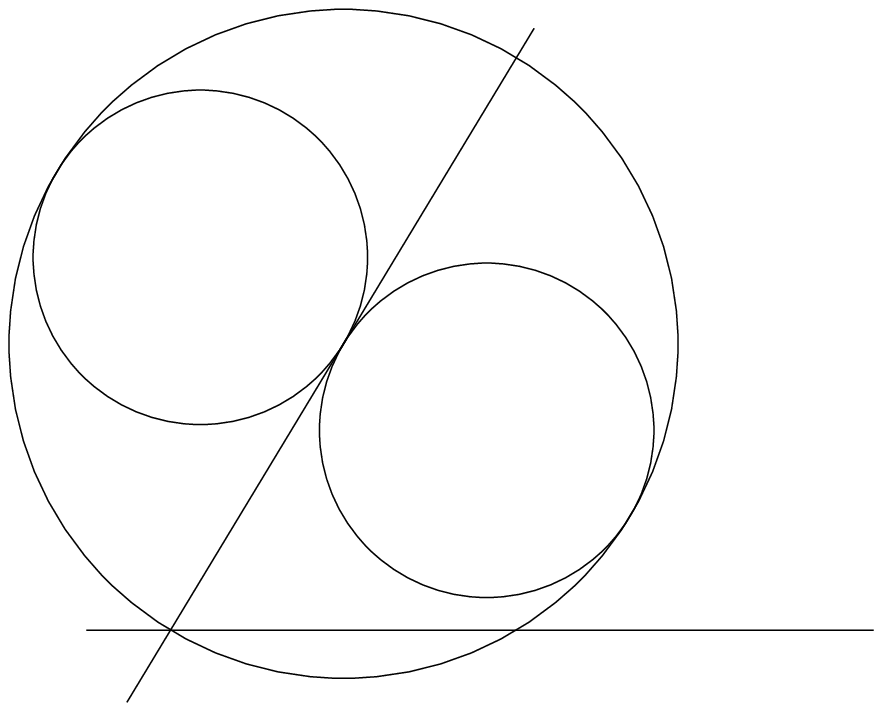}
}
\end{center}
\caption{Critical configuration and beyond}
\end{figure}

It is plain to confirm from Figure \ref{Critical2}
that the case of Figure \ref{TwoCircleMax} appears only when
$\tan(\theta(\X)/2)\le 1/2$. For the case $1/2< \tan(\theta(\X)/2)$,
the maximal ball inscribed in the cell
bounded by $x$-axis and the line of slope $\tan(\theta(\X))$ does not
touch the $x$-axis as in Figure \ref{Larger}.
Therefore the configuration is reduced to the case of
Theorem \ref{2}, i.e., $B(x,r)$ is cut essentially by a single
line, and we get $r/r_n=2$.
All the statements of Theorem \ref{B} are proved.
\end{section}
\bigskip

\begin{center}
{\bf Acknowledgments.}
\end{center}
\bigskip

We would like to express our gratitude to Peter Grabner and
W\"oden Kusner who informed us of the result of Bang, which made Theorem \ref{m+1}
in the present form. The first author is indebted to Hiroyuki Tasaki
for stimulating discussion.
The authors are supported by the
Japanese Society for the Promotion of Science (JSPS), Grant in aid
21540012. The second author expresses his deepest gratitude to the
Hitachi Scholarship Foundation.


\end{document}